\documentclass[11pt]{article}
\usepackage[a4paper,margin=25mm]{geometry}
\usepackage{amsmath,amssymb,amsthm,mathtools,bm}
\usepackage{microtype}
\usepackage{hyperref}
\usepackage{authblk}
\usepackage{fourier}
\usepackage{natbib}

\newtheorem{theorem}{Theorem}

\newtheorem{remark}{Remark}
\theoremstyle{definition}

\newcommand{\gm}{\overline{\varphi}}
\newcommand{\Crit}{\operatorname{Crit}}
\newcommand{\trdeg}{\operatorname{trdeg}}

\title{Unconditional Finiteness of the Number of Critical Points \\ of Gaussian Mixture Densities}
\author[1,2,3]{Akifumi Okuno\thanks{\url{okuno@ism.ac.jp} (corresponding author)}}
\affil[1]{The Institute of Statistical Mathematics}
\affil[2]{The Graduate Institute for Advanced Studies, SOKENDAI}
\affil[3]{RIKEN}
\date{\empty}

\begin{document}
\maketitle

\begin{abstract}
By slightly adapting Wang's recent important argument for the homoscedastic setting, we prove that every finite heteroscedastic Gaussian mixture has only finitely many critical points, and, consequently, only finitely many modes.
\end{abstract}

\section{Introduction}

A wide variety of statistical methods, including clustering and density estimation, use Gaussian mixture densities of the form
\begin{align}
\gm(\bm{x})
=
\sum_{i=1}^k w_i \varphi(\bm{x};\bm{\mu}_i,\bm{\Sigma}_i),
\quad 
\bm{x} \in \mathbb{R}^d, \quad d \in \mathbb{N}, \quad k \in \mathbb{N},
\label{eq:mixture}
\end{align}
where the $i$th Gaussian component is given by
$\varphi(\bm{x};\bm{\mu}_i,\bm{\Sigma}_i) := \det(\bm{\Sigma}_i)^{-1/2} \exp( -(\bm{x}-\bm{\mu}_i)^\top \bm{\Sigma}_i^{-1}(\bm{x}-\bm{\mu}_i)/2 )$.
Here, $w_i>0$, $\bm{\mu}_i\in\mathbb{R}^d$, and $\bm{\Sigma}_i\succ0$ denote arbitrary positive weights, centers, and positive-definite covariance matrices, respectively.

In such applications, each mode is often interpreted as an underlying group, and one might expect the number of modes not to exceed the number of components. However, it is known that the number of modes can exceed the number of components~\citep{perpinan2003isotropic,amendola2019maximum}, and determining the maximum possible number of modes remains an open problem. 
While several works derive upper bounds for general Gaussian mixtures under an a priori finiteness assumption on the modal set~\citep{amendola2019maximum,nguyen2026bounds}, others establish unconditional results under restrictions on either the pair $(d,k)$ or the component geometry~\citep{ray2012upper,edelsbrunner2013isotropic}. In particular, for homoscedastic mixtures with $k=3$, \citet{okuno2026mixtures} went beyond highly symmetric configurations and established unconditional upper bounds. 

Particularly for homoscedastic Gaussian mixture densities, \citet{Wang2026} proved unconditional finiteness of the critical set for arbitrary $k,d \in \mathbb{N}$, although without providing an explicit bound.

\subsection{Main result}
The purpose of this study is to show directly that the homoscedasticity assumption in the argument of \citet{Wang2026} is unnecessary. Our main theorem for heteroscedastic setting is stated as follows.

\begin{theorem}[Heteroscedastic adaptation of \citet{Wang2026}]
\label{thm:main}
For every finite Gaussian mixture density of the form \eqref{eq:mixture} with positive-definite covariance matrices $\bm{\Sigma}_1,\dots,\bm{\Sigma}_k$, the set
\[
\Crit(\gm)=\{\bm{x}\in\mathbb{R}^d:\nabla\gm(\bm{x})=0\}
\]
is finite. Consequently, $\gm$ has finitely many modes.
\end{theorem}

Consequently, the a priori finiteness assumption imposed in \citet{amendola2019maximum} and \citet{nguyen2026bounds} is automatically satisfied. 
See Remark~\ref{rem:ominimal} for why this does not follow from $o$-minimality alone.

\begin{remark}\label{rem:ominimal}
Since $\nabla \gm$ is built from polynomials and the exponential function, $\Crit(\gm)$ is definable in the o-minimal structure $\mathbb{R}_{\exp}$~\citep{wilkie1996model} and hence has finitely many connected components. Thus, Theorem~\ref{thm:main} amounts to excluding positive-dimensional critical components. Neither $o$-minimality nor Khovanskii-type fewnomial and Pfaffian bounds~\citep{khovanskii1991fewnomials} exclude such components by themselves; related Pfaffian approaches were also mentioned in \citet{okuno2026mixtures} but were not used to establish finiteness. This is why \citet{amendola2019maximum} and \citet{nguyen2026bounds} impose an a priori finiteness assumption. A key contribution of \citet{Wang2026} was to overcome this obstruction by introducing an argument based on Ax's functional-transcendence theorem. The present paper simply adapts that argument to the heteroscedastic setting.
\end{remark}

The remainder of this paper is devoted to the proof of Theorem~\ref{thm:main}

\section{Proof of Theorem~\ref{thm:main}}

The proof proceeds by extending the key step in the argument of \citet{Wang2026} from homoscedastic to heteroscedastic Gaussian mixtures. Once this extension is established, the remainder of the proof is identical to that of \citet{Wang2026}.

Suppose that $\Crit(\gm)$ is infinite, toward a contradiction.
For an arbitrary finite Gaussian mixture, \citet{ray2005topography} proved that all critical points are contained in the ridgeline manifold

\[
\mathcal{R} = \left\{
\left(
\sum_{i=1}^{k}\alpha_i\bm{\Sigma}_i^{-1}
\right)^{-1}
\left(
\sum_{i=1}^{k}\alpha_i\bm{\Sigma}_i^{-1}\bm{\mu}_i
\right)
\, \bigg| \, \bm{\alpha}\in\Delta_{k-1} \right\},
\]

where $\Delta_{k-1} = \left\{
\bm{\alpha}\in[0,1]^k \mid 
\sum_{i=1}^k\alpha_i=1 \right\}$ is the probability simplex. Since $\Delta_{k-1}$ is compact and the above parametrization is continuous, $\mathcal{R}$ is compact. Moreover, since $\gm$ is real analytic, its gradient $\nabla\gm$ is continuous, and hence
\[
\Crit(\gm)
=
(\nabla\gm)^{-1}(\{\bm{0}\})
\]
is closed. Therefore, $\Crit(\gm)\subseteq\mathcal{R}$ implies that $\Crit(\gm)$ is compact. Since $\Crit(\gm)$ is infinite by assumption, it contains a non-isolated point $\bm{x}^\ast$.

We use the real-analytic curve selection theorem
\cite[Section~19, Proposition~2]{Lojasiewicz1965}. Since $\Crit(\gm)$ is a real-analytic subset of $\mathbb{R}^d$, there exist $\varepsilon>0$ and a nonconstant real-analytic curve
\[
\bm{\eta}:I\to\Crit(\gm), \qquad I=(0,\varepsilon),
\]
connecting to the non-isolated point $\bm{x}^\ast$ at the endpoint, i.e., $\bm{\eta}(t)\to\bm{x}^\ast$ as $t \searrow 0$.

Here, consider more tractable expression of the mixture density \eqref{eq:mixture}. 
Absorb all factors independent of $\bm{x}$ into positive coefficients $a_i>0$, and set $\bm{A}_i:=\bm{\Sigma}_i^{-1} \succ 0$, $\bm{b}_i:=\bm{A}_i\bm{\mu}_i$, and $g_i(\bm{x}) = -\bm{x}^\top\bm{A}_i\bm{x}/2
+\bm{b}_i^\top\bm{x}$. 
Then, 
\begin{align}
\gm(\bm{x})
=
\sum_{i=1}^k a_i e^{g_i(\bm{x})},
\qquad
\nabla\gm(\bm{x})
=
\sum_{i=1}^k a_i e^{g_i(\bm{x})}
(\bm{b}_i-\bm{A}_i\bm{x}).
\label{eq:grad}
\end{align}

Define $G_i(t):=g_i(\bm{\eta}(t))$ and $Z_i(t):=e^{G_i(t)}$; 
along $\bm{\eta}$, equation $\nabla \gm (\bm{\eta})=\bm{0}$ can be rearranged as $\bigl(\sum_{i=1}^k a_iZ_i\bm{A}_i\bigr)\bm{\eta}=\sum_{i=1}^k a_iZ_i\bm{b}_i$. The matrix on the left is positive definite for every $t\in I$ and is therefore invertible. Thus
\begin{equation}
\bm{\eta}
=
(\eta_1,\eta_2,\ldots,\eta_d)
=
\left(
\sum_{i=1}^k a_i Z_i\bm{A}_i
\right)^{-1}
\left(
\sum_{i=1}^k a_i Z_i\bm{b}_i
\right).
\label{eq:eta}
\end{equation}

For elements $Z_1,\ldots,Z_k$ of a field extension of $\mathbb{R}$, $F:=\mathbb{R}(Z_1,\ldots,Z_k)$ denotes the field generated over $\mathbb{R}$ by $Z_1,\ldots,Z_k$. 
The entries of $\sum_i a_iZ_i\bm{A}_i$ and $\sum_i a_iZ_i\bm{b}_i$ are linear in $Z_1,\dots,Z_k$. Thus the adjugate formula shows that every coordinate of the
right-hand side is a rational function of $Z_1,\dots,Z_k$, and hence $\eta_1,\dots,\eta_d\in F$.

Since each $G_i=-\bm{\eta}^\top\bm{A}_i\bm{\eta}/2+\bm{b}_i^\top\bm{\eta}$ is a quadratic polynomial with real coefficients in the coordinates of $\bm{\eta}$, we also have $G_1,\dots,G_k\in F$. As $e^{G_i}=Z_i$ by definition, adjoining the $G_i$ to $F$ does not enlarge it, and therefore
\begin{align}
F=\mathbb{R}(G_1,\dots,G_k,e^{G_1},\dots,e^{G_k}).
\label{eq:keyfield}
\end{align}

The key difference from \citet{Wang2026} is that equation~\eqref{eq:keyfield} holds not only in the homoscedastic setting considered in \citet{Wang2026}, but also in the heteroscedastic setting considered here. 
Once the identity~\eqref{eq:keyfield} is proved, the remaining argument to prove the contradiction, summarized below for completeness, is identical to that of \citet{Wang2026} Section 2.3. 

Let $v_1,\dots,v_s$ be a maximal subset of $\{G_1,\dots,G_k\}$ that is rationally independent modulo $\mathbb{R}$. By maximality, each $G_i$ is a rational linear combination of $v_1,\dots,v_s$ up to an additive constant, so a positive integral power of each $e^{G_i}$ lies in $F_0:=\mathbb{R}(e^{v_1},\dots,e^{v_s})$; as $F=\mathbb{R}(e^{G_1},\dots,e^{G_k})$, the extension $F/F_0$ is algebraic and $\trdeg_\mathbb{R}F\le s$, where $\trdeg_{\mathbb{R}}$ represents the transcendence degree~\citep[Chapter~VIII]{Lang2002}.  
On the other hand, \eqref{eq:keyfield} yields $\mathbb{R}(v_1,\dots,v_s,e^{v_1},\dots,e^{v_s})\subseteq F$, and Ax's functional-transcendence theorem~\cite[Theorem~3]{Ax1971}, applied in the differential field of real-analytic functions on $I$, bounds the transcendence degree of the left-hand side from below by $s+1$. The two estimates are incompatible, and the theorem follows from the contradiction.
\qed 

\section{Conclusion}
By adapting Wang's recent argument for the homoscedastic setting, we proved that every finite heteroscedastic Gaussian mixture has only finitely many critical points, and hence only finitely many modes.

\section*{Acknowledgement}
A. Okuno was supported by JSPS KAKENHI Grant Numbers 21K17718, 22H05106, and 25K03087. The author thanks Yutaro Kabata for drawing his attention to the recent important work of \citet{Wang2026}.

\section*{Declaration of AI use}
AI assistance was used for mathematical discussion and language editing. In particular, GPT-5.6 Sol Pro was used to summarize the work of \citet{Wang2026}, isolate the part of its argument essential to the present problem, and develop a streamlined route to the proof. Subsequent interactions with the model were used to refine the central observation that the identity in \eqref{eq:keyfield} extends from the homoscedastic to the heteroscedastic setting. The author subsequently verified every mathematical claim and argument. 

\bibliographystyle{apalike}
\bibliography{references}

\end{document}